\documentclass[12pt]{article}
\usepackage{amsmath,amscd,amsthm,amsfonts,amssymb}
\usepackage{lipsum}
\usepackage{graphics}
\usepackage{amsmath,amscd,amsthm,amsfonts,amssymb} 
\usepackage{pgfpages}
\usepackage[shortlabels]{enumitem}
\usepackage{adjustbox}
\usepackage{graphicx}
\usepackage{tikz}
\usepackage{caption}\usepackage{mdwlist}
\usepackage{mhchem}
\begin{document}
	\newtheorem{notat}{Notation}
	\newtheorem{fact}{Fact}
	\newtheorem{ass}{Assumption}
	\newtheorem{thm}{Theorem}
	\newtheorem{cor}{Corollary}
	\newtheorem{lem}{Lemma}
	\newtheorem{slem}{Sublemma}
	\newtheorem{prop}{Proposition}
	\newtheorem{defn}{Definition}
	\newtheorem{conj}{Conjecture}
	\newtheorem*{ques*}{Question}
	\newtheorem{claim}{Claim}
	\newcounter{constant} 
	\newtheorem{rmk}{Remark}
	\newtheorem*{rmk*}{Remark}
	\newcommand{\newconstant}[1]{\refstepcounter{constant}\label{#1}}
	\newcommand{\useconstant}[1]{\Delta_{\ref{#1}}}
	\title{Casson invariant of knots defined by random braids}
	\date{}
	\author{Marina Ville}
\maketitle
\begin{abstract}
We investigate the expectation $E(v_2)$ of the Casson invariant $v_2$ for knots represented by two similar types of braids, the rosette $s$-braids which are products of $k$ factors $\prod_{1\leq 2i\leq s}\sigma_{2i}^{\pm 1}\prod_{1\leq 2i+\leq s}\sigma_{2i+1}^{\pm 1}$
and the  quasitoric $s$-braids which are products of $k$ factors $\sigma_{1}^{\pm 1}\sigma_{2}^{\pm 1}...\sigma_{s-1}^{\pm 1}$.
If $k$ is of the form $ns+1$, $\dfrac{E(v_2)}{c}\sim \dfrac{s+1}{24}$ where $c$ is the number of crossing points of  braid; if $k$ is of the form $ns-1$, $\dfrac{E(v_2)}{c}\sim \dfrac{5-s}{24}$.\\
If we go through all $3$-braids, regardless of the shape of their diagram, we get $\dfrac{E(v_2)}{c}\sim \dfrac{1}{8}$.   \\
For knots with a checkerboard diagram, $\dfrac{E(v_2)}{c}\sim \dfrac{1}{16}$. 
\end{abstract}
\section{Introduction}
\subsection{Overview}
Knots form a discrete  set so we would like to do statistics on them. Since there are infinitely many knot types, one  divides the knot types in finite classes $C_n$ (indexed by integers $n$) of equiprobable knots. To study a given knot invariant $X$, we compute its average $E(X,C_n)$ on each $C_n$ and see how $E(X,C_n)$ behaves as $n$ tends to infinity. For example, if we fix a knot diagram $\Gamma_n$ for every $n$, $C_n$ will be the set of knots projecting on the plane to $\Gamma_n$: for each crossing point of $\Gamma_n$, a given strand of a knot in $C_n$ will be above or under the other one with probability $\dfrac{1}{2}$. \\
\\ We do this with braid shadows, the shadow of a braid $\beta$ being its projection on the horizontal plane; equivalently it is the data of the braid generators $\sigma_i$'s appearing in $\beta$ without the exponents $\pm 1$. We look at two types of braids, {\it rosette braids} and {\it quasitoric braids}. These are very similar braids which  both represent torus knots; the literature says that they are conjugate but we have not found a general proof so we study both of them and get identical results.   \\
Let $s$ be a positive integer. Lamm ([La]) defined  a {\it rosette $s$-braid} of type $(s,k)$, denoted ${\mathcal R}(s,k)$, as a product of $k$ factors of the form  $$\prod_{1\leq 2i\leq s}\sigma_{2i}^{\pm 1}\prod_{1\leq 2i+\leq s}\sigma_{2i+1}^{\pm 1}.$$
Below is the closure of a $\mathcal R(4,5)$ braid diagram.\\
\\	
\begin{tikzpicture}[thick, scale=0.5]
\draw [very thick] (0,1) -- (-0.5,1);
\draw [very thick] (0,2) -- (-1,2);
\draw [very thick] (0,3) -- (-1.5,3);
\draw [very thick] (0,4) -- (-2,4);

\draw [very thick] (10,1) -- (10.5,1);
\draw [very thick] (10,2) -- (11,2);
\draw [very thick] (10,3) -- (11.5,3);
\draw [very thick] (10,4) -- (12,4);

\draw [very thick] (-0.5,0) -- (-0.5,1);
\draw [very thick] (-1,-1) -- (-1,2);
\draw [very thick] (-1.5,-2) -- (-1.5,3);
\draw [very thick] (-2,-2.5) -- (-2,4);

\draw [very thick] (10.5,0) -- (-0.5,0);
\draw [very thick] (11,-1) -- (11,2);
\draw [very thick] (11.5,-2) -- (11.5,3);
\draw [very thick] (12,-2.5) -- (12,4);

\draw [very thick] (10.5,0) -- (10.5,1);
\draw [very thick] (11,-1) -- (-1,-1);
\draw [very thick] (11.5,-2) -- (-1.5,-2);
\draw [very thick] (12,-2.5) -- (-2,-2.5);

\draw [very thick] (0,1) -- (3,4);
\draw [very thick] (0,3) -- (1,4);
\draw [very thick] (0,4) -- (3,1);
\draw [very thick] (0,2) -- (1,1);
\draw [very thick] (2,4) -- (1,4);
\draw [very thick] (2,1) -- (1,1);

\draw [very thick] (2,1) -- (5,4);
\draw [very thick] (2,3) -- (3,4);
\draw [very thick] (2,4) -- (5,1);
\draw [very thick] (2,2) -- (3,1);
\draw [very thick] (4,4) -- (3,4);
\draw [very thick] (4,1) -- (3,1);

\draw [very thick] (4,1) -- (6,3);
\draw [very thick] (4,3) -- (5,4);
\draw [very thick] (4,4) -- (6,2);
\draw [very thick] (4,2) -- (5,1);
\draw [very thick] (6,4) -- (5,4);
\draw [very thick] (6,1) -- (5,1);

\draw [very thick] (6,1) -- (8,3);
\draw [very thick] (6,3) -- (7,4);
\draw [very thick] (6,4) -- (9,1);
\draw [very thick] (6,2) -- (7,1);
\draw [very thick] (8,4) -- (7,4);
\draw [very thick] (8,1) -- (7,1);

\draw [very thick] (8,1) -- (10,3);
\draw [very thick] (8,3) -- (9,4);
\draw [very thick] (8,4) -- (10,2);
\draw [very thick] (8,2) -- (9,1);
\draw [very thick] (10,4) -- (9,4);
\draw [very thick] (10,1) -- (9,1);

\end{tikzpicture}	\\

Manturov ([Ma]) defined a {\it quasitoric} $s$-braid of type $(s,k)$, denoted ${\mathcal Q}(s,k)$ as a product of $k$ factors of the form  
 $$\sigma_{1}^{\pm 1}\sigma_{2}^{\pm 1}...\sigma_{s-1}^{\pm 1}.$$
Below is the closure of a $\mathcal Q(4,5)$ braid diagram.\\

 \begin{tikzpicture}[thick, scale=0.5]
 \draw [very thick] (0,1) -- (-0.5,1);
 \draw [very thick] (0,2) -- (-1,2);
 \draw [very thick] (0,3) -- (-1.5,3);
 \draw [very thick] (0,4) -- (-2,4);
 
 \draw [very thick] (10,1) -- (10.5,1);
 \draw [very thick] (10,2) -- (11,2);
 \draw [very thick] (10,3) -- (11.5,3);
 \draw [very thick] (10,4) -- (12,4);
 \draw [very thick] (0,1) -- (6,4);
 \draw [very thick] (0,2) -- (4,4);
 \draw [very thick] (0,3) -- (2,4);
 
 \draw [very thick] (2,1) -- (8,4);
 
 \draw [very thick] (4,1) -- (10,4);
 
 \draw [very thick] (0,4) -- (2,1);
 \draw [very thick] (2,4) -- (4,1);
 \draw [very thick] (4,4) -- (6,1);

 \draw [very thick] (6,4) -- (8,1);
 \draw [very thick] (8,4) -- (10,1);
 \draw [very thick] (6,1) -- (10,3);
 \draw [very thick] (8,1) -- (10,2);

 \draw [very thick] (10,1) -- (10.5,1);
 \draw [very thick] (10,2) -- (11,2);
 \draw [very thick] (10,3) -- (11.5,3);
 \draw [very thick] (10,4) -- (12,4);
 
 \draw [very thick] (-0.5,0) -- (-0.5,1);
 \draw [very thick] (-1,-1) -- (-1,2);
 \draw [very thick] (-1.5,-2) -- (-1.5,3);
 \draw [very thick] (-2,-2.5) -- (-2,4);
 
 \draw [very thick] (10.5,0) -- (-0.5,0);
 \draw [very thick] (11,-1) -- (11,2);
 \draw [very thick] (11.5,-2) -- (11.5,3);
 \draw [very thick] (12,-2.5) -- (12,4);

 \draw [very thick] (10.5,0) -- (10.5,1);
 \draw [very thick] (11,-1) -- (-1,-1);
 \draw [very thick] (11.5,-2) -- (-1.5,-2);
 \draw [very thick] (12,-2.5) -- (-2,-2.5);

 \end{tikzpicture}
 \\
The plat closure of a rosette braid projects onto a checkerboard diagram. A $(p,q)$ checkerboard diagram is parametrized by a Lissajous curve\\ $x=\cos(qt+\alpha), y=\cos(pt)$ with $(p,q)=1$. Several authors have proved that checkerboards diagrams are models for random knots ([EZ]). Below is a $(3,7)$ checkerboard diagram.

	\begin{tikzpicture}
\draw [very thick] (0,0) -- (2.5,2.5);
\draw [very thick] (1,0) -- (3.5,2.5);
\draw [very thick, -stealth] (1,0) -- (1.5,0.5);
\draw [very thick] (-0.5,0.5) -- (-0.5,0);
\draw [very thick] (0,0) -- (-0.5,0);
\draw [very thick] (-0.5,1) -- (0.5,0);
\draw [very thick] (-0.5,0.5) -- (1.5,2.5);
\draw [very thick] (-0.5,1.5) -- (-0.5,1);
\draw [very thick] (0,1.5) -- (1.5,0);
\draw [very thick] (-0.5,1.5) -- (0.5,2.5);
\draw [very thick] (1,2.5) -- (0.5,2.5);
\draw [very thick] (1,0) -- (0.5,0);
\draw [very thick,-stealth] (3.5,0)--(2,1.5);
\draw [very thick] (3.5,0)--(1,2.5);
\draw [very thick, -stealth] (0,2.5) -- (2.5,0);
\draw [very thick] (2,0) -- (4.5,2.5);
\draw [very thick] (2,0) -- (1.5,0);
\draw [very thick] (2.5,0) -- (3,0);
\draw [very thick] (3.5,2.5) -- (4,2.5);
\draw [very thick] (2,2.5) -- (1.5,2.5);
\draw [very thick] (2.5,2.5) -- (3,2.5);
\draw [very thick] (3.5,2.5) -- (4,2.5);
\draw [very thick, -stealth] (2,2.5) -- (4.5,0);
\draw [very thick] (3,2.5) -- (5.5,0);
\draw [very thick] (4,2.5) -- (6,0.5);
\draw [very thick] (4.5,0) -- (5,0);
\draw [very thick,-stealth] (4,0) -- (3.5,0);
\draw [very thick] (0,1) -- (1.5,2.5);
\draw [very thick] (3,0) -- (5.5,2.5);
\draw [very thick] (0,2) -- (0.5,2.5);
\draw [very thick] (4.5,2.5) -- (5,2.5);
\draw [very thick] (6,1.5) -- (5,2.5);
\draw [very thick] (4.5,2.5) -- (5,2.5);
\draw [very thick] (4,0) -- (6,2);
\draw [very thick,-stealth] (6,2) --(5,1);
\draw [very thick, -stealth] (5.5,2.5) -- (6,2.5);
\draw [very thick] (5,0) -- (6,1);
\draw [very thick] (6,0.5) -- (6,0);
\draw [very thick, -stealth] (6,2.5) -- (6,2);
\draw [very thick] (5.5,0) -- (6,0);
\draw [very thick, -stealth] (-0.5,2.5) -- (0,2.5);
\draw [very thick] (-0.5,2.5) -- (-0.5,2);
\draw [very thick] (1.5,0) -- (-0.5,2);

\draw [very thick] (6,1) -- (6,1.5);
\end{tikzpicture}\\
\\
We estimate the expectancy $E(v_2)$ of the Vassiliev invariant $v_2$ for these three models. A Vassiliev invariant of knots/links extends to singular knots/links and is of order $p$ if it vanishes for links with more than $p$ singular points. The simplest one is the linking number. Next comes $v_2$ which is of order $2$ and coincides with the coefficient of $z^2$ in the Conway polynomial. It is also called the Casson invariant of knots because it appears in the theory of the Casson invariant of homology $3$-spheres.\\ For example,
 [EZ-H-L-N] studies the {\it petaluma model} where the knot diagram looks like a daisy with $2n+1$ petals; for this model, $E(v_2)=\dfrac{n^2}{24}\sim \dfrac{c}{48}$. Here and in the rest of the paper, $c$ denotes the number of crossing points of a diagram; for the petaluma diagram, it is the number of crossing points of a desingularization of the diagram.\\
 We let $\hat{\sigma}_i$ be the shadow  of a braid generator $\sigma_i$: when we lift the braid shadow to a regular braid,  we give a sign $+1$ or $-1$ to each $\hat{\sigma_i}$ with equal probability.\\
 \\
 Following Lamm and Manturov (\S \ref{definition des rosettes}), every knot type can be described by a rosette braid and a quasitoric braid of type $(s,as+1)$  where $a\in\mathbb{N}$; for these $(s,as+1)$ models  $\dfrac{E(v_2)}{c}\sim\dfrac{s+1}{24}$ (proved by [C-E1] in the quasitoric case).  Similarly, a knot can be represented by a  rosette braid and a quasitoric braid of type $(s,sb-1)$, with $b\in\mathbb{N}$; for both these models, we have  $\dfrac{E(v_2)}{c}\sim\dfrac{5-s}{24}$. Thus the $(s,sn+1)$ and $(s,sn-1)$ models are close but they give very different expectations. Note that [C-E1] and Chaim Even Zohar exhibited quasitoric braids with negative $E(v_2)$.\\
 \\
 Then we look at $3$-braids since their computations are more tractable.
The expectations $E(v_2)$ for $(3,3a+1)$ and $(3,3a-1)$ are respectively equivalent to $\dfrac{c}{6}$ and $\dfrac{c}{12}$. If we now go through all the $3$-braids, with no restriction on their shape, we get $\dfrac{E(v_2)}{c}\sim\dfrac{1}{8}$ which is the average of $\dfrac{1}{6}$ and $\dfrac{1}{12}$.\\
\\
Finally we consider knots projecting onto a checkerboard diagram and compute $\dfrac{E(v_2)}{c}\sim \dfrac{1}{16}$. 

\subsection{The Casson invariant}

\subsubsection{A combinatorial description}

Polyak and Viro ([P-V]) described $v_2$ in terms of {\it Gauss diagrams} where the knot is encoded by its crossing points written along a circle in the order that they occur when we travel along the knot. The two points on the knot corresponding to a given crossing point are joined by a chord; an arrow on this chord encodes which strand is above and which strands is below. Polyak and Viro derive a  combinatorial way of computing $v_2$ which we recall now.\\ 
Let $K$ be an oriented knot and fix a starting point on $K$. A crossing point of $K$ is {\it ascending} (resp. {\it descending}) if the parameter goes through the lower strand  before (resp. after) going through its upper strand. \\
NOTATION. If $(p,\tilde{p})$ is a double point of $K$, $p$ (resp. $\tilde{p}$) denotes the first (resp. second) point of $K$ which is encountered by the parameter. 
\begin{prop}\label{definition de cdeux}
	A pair of crossing points $((p,\tilde{p}),(q,\tilde{q}))$ of an oriented knot is said to be valid if 
	\begin{enumerate}
		\item $((p,\tilde{p})$ is ascending
		\item $(q,\tilde{q})$ is descending
		\item The crossing points are interleaved, i.e. they  appear in the order $p,q,\tilde{p},\tilde{q}$. 
	\end{enumerate}
	Then $$v_2(K)=\sum sign(p)sign(q)$$ where $((p,\tilde{p}),(q,\tilde{q}))$ goes through all valid pairs of crossing points.
\end{prop} 
Recall that the signs of crossing points of an oriented knot are
\begin{equation}\label{signe des points doubles}
\begin{tikzpicture}
\draw[very thick](0,1)--(0.45,0.55);
\draw[very thick,-stealth](0.55,0.45)--(1,0);
\draw[very thick,-stealth](0,0)--(1,1);
\end{tikzpicture}-1.
\ \ \ \ \ \ \ \ \ \ \ \ \ \ \ \
\begin{tikzpicture}
\draw[very thick,-stealth](0,1)--(1,0);
\draw[very thick](0,0)--(0.45,0.45);
\draw[very thick,-stealth](0.55,0.55)--(1,1);
\end{tikzpicture} +1.
\end{equation}

Note that [EZ-H-L-N] derive Proposition \ref{definition de cdeux} from the characterization of $v_2$ as the coefficient of $z^2$ in the Conway polynomial, using the skein relations for the Conway polynomial.\\ 
REMARK. In Proposition \ref{definition de cdeux} above, 1. and 2. depend on the sign of the crossing points but 3. only depends on their position on the braid shadow. \\
To derive the expectation of $v_2$, we define 

\begin{defn}\label{epsilon} (cf. [B-B]) Let $\Gamma$ be a closed plane curve with transverse double points and fix a starting point  $x\in\Gamma$. 
If $d$ is a crossing point of $\Gamma$, the number 	$\epsilon_d(x)$ is $+1$ if the outgoing branches of
$\Gamma$ at $d$ (in the order in which they are visited when travelling along the curve starting at $x$) orient the plane clockwise; $\epsilon_d(x)$ is $-1$ otherwise.\\
	Let $\hat{\Gamma}$ be a knot projecting on $\Gamma$. A double point $(d,\tilde{d})$ of $\hat{\Gamma}$ is ascending (resp. descending) if $sign(d)=\epsilon_d(x)$ (resp. $sign(d)=-\epsilon_d(x)$)
\end{defn}
Here is an illustration with two ascending points; $\textcircled{1}$ denotes the first strand travelled by the parameter and $\textcircled{2}$ denotes the second one.\\
\begin{tikzpicture}
\draw[very thick](0,1)--(0.45,0.55);
\draw[very thick,-stealth](0.55,0.45)--(1,0);
\draw[very thick,-stealth](0,0)--(1,1);
	\node at (-0.5,1){\textcircled{1}};
\node at (-0.5,0){\textcircled{2}};
\node at (0.5,-0.5){$\epsilon_d(x)=sign(d)=-1$};
\end{tikzpicture}\ \ \ \ \ \ \ \ \ \ \ \ \
\begin{tikzpicture}
\draw[very thick,-stealth](0,1)--(1,0);
\draw[very thick](0,0)--(0.45,0.45);
\draw[very thick,-stealth](0.55,0.55)--(1,1);
	\node at (-0.5,1){\textcircled{2}};
\node at (-0.5,0){\textcircled{1}};
\node at (0.5,-0.5){$\epsilon_d(x)=sign(d)=1$};
\end{tikzpicture} \\

We consider a knot diagram $\Gamma$ and order its crossing points by denoting them $p_1,...,p_c$. We let ${\mathcal K}(\Gamma)$ be the set of knots projecting on $\Gamma$.\\
The cardinal number of ${\mathcal K}(\Gamma)$ is $2^c$; its elements all have the same probability, and we compute the expectation $E(v_2)$ of the invariant $v_2$.\\
Note that interleaved crossing points $((p,\tilde{p})$, $(q,\tilde{q}))$ form a valid pair if and only if $sign(p)=\epsilon_p(x)$ and  $sign(q)=-\epsilon_q(x)$. We derive
\begin{prop}\label{weighted sums}
	$E(v_2)=\dfrac{1}{2^c}\sum_{K\in {\mathcal K}(\Gamma)} v_2(K)=-\dfrac{1}{4}\sum_{p,q}\epsilon_p(x)\epsilon_q(x)$ where $p,q$ go through the interleaved pairs of crossing points of $\Gamma$. 
\end{prop}
\subsection{Connection with the Arnold invariants for closed plane curves}\label{arnold invariants}

	The expectation of the Casson invariants for a given knot diagram $\Gamma$ can be expressed in terms of the Arnold invariants for closed plane curves in $\mathbb{R}^2$ with transverse double points.
\begin{equation}
E(v_2)=\dfrac{1}{8}defect(\Gamma)
\end{equation}
where the {\it defect} is the combination of two other Arnold invariants, $J^+$ and $St$,
\begin{equation}
defect=J^++2St
\end{equation}
Using [B-B-P], we recall these invariants (see also some of the original papers [E-C], [E-C-C], [Lu], [Po], [Sh]).\\
The {\it index} of a region $U$ in $\mathbb{R}^2-\Gamma$ is the winding number of $\Gamma$ around a point of $U$. If $p$ is a point in $\Gamma$, $ind_\Gamma(p)$ is the average of the indices of the $2$ regions of $\mathbb{R}^2-\Gamma$ adjacents to $p$ ($4$ regions if $p$ is a double point). If we fix a starting point $x$  of $\Gamma$, we define
\begin{equation}
St=\sum_{d\ \mbox{double point of}\ \Gamma}\epsilon_d(x)ind_K(d)+ind_K(x)^2-\dfrac{1}{4}
\end{equation}
We now consider a Seifert circle $c$ of $\Gamma$ and we let $index(c)$ be $+1$ if the circle is
oriented counterclockwise and $-1$ otherwise. We define $t(c)$ as
\begin{equation} 
t(c) = index(c)ind_\Gamma(R_c)
\end{equation}
where $R_c$  is a region  which borders some arc of $c$ from the
outside. Then
\begin{equation}
J^+=1+\nu-s-2\sum_{\mbox{c Seifert circle of}\ \Gamma}t(c)
\end{equation}
where $\nu$ (resp. $s$) is the number of crossing points (resp. Seifert circles) of $\Gamma$.
\subsection{Models of knots}\label{definition des rosettes}
Lamm and Manturov showed that knots and links can be expressed as closures of rosette and quasitoric braids. 
\begin{thm}([La], [Ma]) Let $\alpha\in B_s$ be a braid which closes in a knot. It can be represented by the following $4$ braids: rosette braids $\mathcal R(s,sa+1)$, $\mathcal R(s,sb-1)$ and quasitoric braids $\mathcal Q(s,sc+1)$, $\mathcal Q(s,sd-1)$ where $a,b,c,d\in\mathbb{N}$. 
	 \end{thm}
 \begin{proof} The statement is slightly different from Lamm and Manturov's statement so we tweak their proof a bit to obtain it but the ingredients are in their papers.
 	\begin{notat}\label{permutations} We number the $s$ endpoints of a $s$-braid $\beta$ from top ($1$-st strand) to bottom ($s$-th strand). We define a permutation $\pi(\beta)\in \mathfrak{S}_s$ by setting: the strand starting at the $i$-th point ends at the $\pi(i)$-th point.
 	\end{notat}
  Recall that a braid $\gamma$ is {\it pure} if $\pi(\gamma)=Id$. Lamm and Manturov prove that every pure braid can be represented by a rosette braid $\mathcal R(s,sn)$ and a quasitoric braid $\mathcal Q(s,sm)$ for some integers $n,m$.\\
Now let $\alpha$ be a $s$-braid which  closes in a knot; equivalently $\pi$ is a $1$-cycle. Let $\delta$ be a braid (rosette or quasitoric) of type $(s,1)$ (resp. $(s,s-1)$). 	Both $\pi(\alpha)$ and $\pi(\delta)$ are $1$-cycles, hence they are conjugate; thus for some braid $\beta$, we have 
\begin{equation}
\pi(\alpha)=\pi(\beta)\pi(\delta)\pi(\beta)^{-1}
\end{equation}
\begin{equation}
\mbox{hence}\ \ \ \ \ \pi(\delta)=\pi(\beta)^{-1}\pi(\alpha)\pi(\beta) \ \ \ \ \ \ \ \ \ \ \ \
\end{equation}
Thus $\pi(\delta^{-1}\beta^{-1}\alpha\beta))=Id$ so $\delta^{-1}\beta^{-1}\alpha\beta$ is a pure braid and  it is represented by a rosette/quasitoric braid $Z$ of type $(s,ns)$ for some $n$, hence $\beta\alpha\beta^{-1}$ is represented by $\delta Z$ which is a braid of type $(s,ns+1)$ (resp. $(s,ns+s-1)$). Thus the knot represented by $\alpha$ is represented both by a braid 
$(s,ns+1)$ and by a braid $(s,ns+s-1)$ (in other words by a braid of type $(s,ms-1)$).
\end{proof}

\subsection{The results}
We compute the expectation of the Casson invariant for the knots projecting on the rosette/quasitoric  braids of type $\mathcal R(s,ns+1)$, $\mathcal Q(s,ns+1)$ and $\mathcal R(s,ns-1)$, $\mathcal Q(s,ns-1)$. The former is always positive and the latter is negative for $s>5$. The case $\mathcal Q(s,ns+1)$ has been proved by [C-E1] using Arnold invariants but we reprove it to use it for the case $\mathcal Q(s,ns-1)$.

\begin{thm}\label{thm fermeture usuelle}
Given $s,n$, the expectation of the Casson invariant over the knots given by the closure of a rosette or quasitoric braid of type $(s,ns+1)$ (resp. $(s,ns-1)$) is given by 
$$  E(v_2)=\dfrac{s(s^2-1)n}{24}\ \ \ \ \ \ \ \ \ \Big(\mbox{resp.}\ \  E(v_2)=\dfrac{s(s-1)(5-s)n}{24} \Big)$$
\end{thm}
 If $c$ is the length of the rosette braid $\mathcal R(s,ns+1)$, $\mathcal Q(s,ns+1)$ (resp. $\mathcal R(s,ns-1)$, $\mathcal Q(s,ns-1)$), we have $\lim\limits_{c\longrightarrow+\infty}\dfrac{E(v_2)}{c}=\dfrac{s+1}{24}$ (resp. $\lim\limits_{c\longrightarrow+\infty}\dfrac{E(v_2)}{c}=\dfrac{5-s}{24}$).\\
\\
We can view $(s,ns+1)$-braids and $(s,ns-1)$-braids as extreme cases for $s$-braids. For $3$-braids, they respectively yield 
$E(v_2)\sim\dfrac{c}{6}$ and $E(v_2)\sim\dfrac{c}{12}$. Now, instead of limiting ourselves to rosette or quasitoric braid shadows, we compute the expectation of the Casson invariant over all $3$-braid shadows.
\begin{prop}\label{autre modele}
	Let $c$ be an even number. Consider the set of equiprobable braid shadows $\gamma$ of the same length $c$ which verify
	\begin{itemize}
		\item
		$\gamma=\hat{\sigma}_{a_1}...\hat{\sigma}_{a_c}$ where each $a_{i_k}$ belongs to $\{1,2\}$ 
		\item
		$\gamma$ closes in a knot
	\end{itemize}
	For each braid shadow above, we give a sign $\pm 1$ to each $\sigma_{a_1}$ with equal probability $\dfrac{1}{2}$. Then 
	\begin{equation}\label{asymptote1} \lim\limits_{c\longrightarrow+\infty}\dfrac{E(v_2)}{c}=\dfrac{1}{8}. 
	\end{equation}

\end{prop}
 The following theorem comes from [K-O] but the proof has been lost so we reprove it. 
\begin{thm}\label{thm fermeture plat}
	The expectation of the Casson invariant over the knots given by a $(p,q)$ checkerboard diagram is
 $$E(v_2)=\dfrac{(p-1)(q-1)}{8}$$
\end{thm}

\section*{Acknowledgements} I am very grateful to Tahl Nowik who suggested this line of research, to Chaim Even Zohar and Moshe Cohen for very helpful conversations and to Toru Ohmoto for kindly answering my question about his work. I  thank both Judy Kupferman and ICTS, Bengaluru for their hospitality during part of this work.  
\section{Proof of Theorem \ref{thm fermeture usuelle}}
\subsection{The braid shadow}\label{notations usual closure}
The building block of the rosette $s$-braid shadows is $X=\prod \hat{\sigma}_{2i+1}\prod \hat{\sigma}_{2i}$,\\ for example; for $6$ strands, $X$ is\ \ \ \ \
\begin{minipage}{0.2\textwidth}
	\begin{tikzpicture}
	\draw[very thick](0,2.5)--(1,1.5); 
	\draw[very thick](0.5,2.5)--(1,2.5);
	\draw[very thick](0.5,2.5)--(0,2);
	\draw[very thick](0,1.5)--(1,0.5);
	\draw[very thick](0,1)--(1,2);
	\draw[very thick](0,0)--(1,1);
	\draw[very thick](0,0.5)--(0.5,0);
	\draw[very thick](1,0)--(0.5,0);
	\end{tikzpicture}
\end{minipage}
\begin{minipage}{0.5\textwidth}
Note that	$X^s$ is a pure braid; we consider here $X^{ns+1}=(X^s)^nX$ and $X^{ns-1}=(X^s)^{n-1}X^{s-1}$.
\end{minipage}
For quasitoric braids, the building block $Y$ is
\begin{minipage}{0.1\textwidth}
	\begin{tikzpicture}
	\draw[very thick](0,2.5)--(1,0); 
	\draw[very thick](0,2)--(1,2.5);
	\draw[very thick](1,2)--(0,1.5);
	\draw[very thick](1,1.5)--(0,1);
	
	\draw[very thick](0,0.5)--(1,1);
	\draw[very thick](1,2)--(0,1.5);
	\draw[very thick](0,0)--(1,0.5);
	\end{tikzpicture}

\end{minipage}
\begin{minipage}{0.6\textwidth}
Similarly to above, $Y^s$ is a pure braid\\ and we consider $Y^{ns+1}$ and  $Y^{ns-1}$
\end{minipage}

\subsection{Permutations}
\subsubsection{General case}\label{ou est le starting point}
We consider a $s$-braid $\beta$ and recall the permutation $\pi$ defined in Notation \ref{permutations}.
 We parametrize $\beta$ by a parameter $t:[0,1]\longrightarrow \mathbb{R}^2$ starting at the starting point of the first strand (left up). We define $s$ intervals $I_k$, $1\leq k\leq s$,  such that\\
1) $t$ starts at the beginning of $I_1$, goes through $I_{k-1}$ just before $I_k$ \\ 2) each $t(I_k)$ is a strand of the closed braid and we say that it is the $\phi(k)$-th strand and define the permutation $\phi$ of $\{1,...,s\}$.
\begin{lem}\label{lien entre les deux}
If $i\in\{1,...,s\}$,	$\phi(i)=\pi^{i-1}(1)$
\end{lem}
\begin{proof}
We prove it by induction. We know that $\phi(1)=1=\pi^{0}(1)$. If $\phi(i)=\pi^{i-1}(1)$, the $i$-th interval starts at the point $\pi^{i-1}(1)$, so it ends at $\pi^{i}(1)$ so the $(i+1)$-th interval starts at $\pi^{i}(1)$ and $\phi(i+1)=\pi^{i}(1)$.
\end{proof}

\subsubsection{The case $\mathcal R(s,1)$}
We let $\phi_1$ be the permutation $\phi$ defined above for a braid in  $\mathcal R(s,sn+1)$, equivalently for $\mathcal R(s,1)$.
\begin{lem}\label{forme de u} Let $A=\Big[\dfrac{s+1}{2}\Big]$. Then
	\begin{enumerate}
		\item $\phi_1(1)=1$
		\item if  $i\leq A$, then $\phi_1(i)=2i-1$
		\item if  $i\geq A+1$, then $\phi_1(i)=2s+2-2i$
	\end{enumerate}
\end{lem}

\begin{proof}. The figures show how things work.\\
	
	\begin{minipage}{2.7cm}
		\begin{tikzpicture}
		\node at (-0.5,2.5) {$\phi_1(1)$};
		\draw[very thick](0,2.5)--(1,1.5); 
		\node at (-0.5,2) {$\phi_1(6)$};
		
		\node at (-0.5,1.5) {$\phi_1(2)$};
		\draw[very thick](0.5,2.5)--(1,2.5);
		
		\draw[very thick](0.5,2.5)--(0,2);
		\node at (-0.5,0.5) {$\phi_1(3)$};
		\draw[very thick](0,1.5)--(1,0.5);
		
		\draw[very thick](0,1)--(1,2);
		
		\draw[very thick](0,0)--(1,1);
		\node at (-0.5,1) {$\phi_1(5)$};
		\draw[very thick](0,0.5)--(0.5,0);
		\node at (-0.5,0) {$\phi_1(4)$};
		\draw[very thick](1,0)--(0.5,0);
		\end{tikzpicture}\\
		$s$ is even 
	\end{minipage}
	\begin{minipage}{2.7cm}
		\begin{tikzpicture}
		\node at (-0.5,2.5) {$\phi_1(1)$};
		\draw[very thick](0,2.5)--(1,1.5); 
		\node at (-0.5,2) {$\phi_1(5)$};
		
		\node at (-0.5,1.5) {$\phi_1(2)$};
		\draw[very thick](0.5,2.5)--(1,2.5);
		
		\draw[very thick](0.5,2.5)--(0,2);
		\node at (-0.5,0.5) {$\phi_1(3)$};
		\draw[very thick](0,1.5)--(1,0.5);
		
		\draw[very thick](0,1)--(1,2);
		
		\draw[very thick](0.5,0.5)--(1,1);
		\node at (-0.5,1) {$\phi_1(4)$};
		\draw[very thick](0.5,0.5)--(0,0.5);

		\end{tikzpicture}\\
		$s$ is odd
	\end{minipage}
	\begin{minipage}{7.5cm}
		We have $\phi_1(1)=1$. After that, every time we go from the beginning of the interval $I_i$ to its end (which is also the beginning of the interval $I_{i+1}$), we go down $2$ strands. Thus $\phi_1(i)=\phi_1(i-1)+2$ and by induction $\phi_1(i)=2i-1$ (2. of the Lemma). This works until (included) $2i-1=s$ or $2i=s$, i.e. $i=A$. If $s$ is odd (resp. even), $\phi_1(A)=s$,  $\phi_1(A+1)=s-1$ (resp. $\phi_1(A)=s-1$,  $\phi_1(A+1)=s
		$).\\
	 Then we go upwards by $2$ strands every time, we derive $\phi_1(i+1)=\phi_1(i)-2$ and we check 3. above.
	\end{minipage}
	
\end{proof}
\subsubsection{The case $\mathcal R(s,s-1)$}\label{s,s moins un}
We denote by $\pi_{-1}$, $\phi_{-1}$ the permutations defined by Notation \ref{permutations} for $\mathcal R(s,ns-1)$, equivalently $\mathcal R(s,s-1)$. We have $\pi_{-1}=\pi_1^{s-1}$ and since $(s,s)$ is a pure braid, $\pi_1^{s}=Id$ and $\pi_{-1}=\pi_1^{-1}$. Thus\\
$\bullet$ $\phi_{-1}(1)=1$, \\
$\bullet$ if $i\geq 2$, $\phi_{-1}(i)=\pi_{-1}^{i-1}(1)=\pi_1^{1-i}(1)=\pi_1^{s+1-i}(1)= \phi_1(s+2-i).$\\
\\
We illustrate this by the following image of the pure braid $X^s$, where the non dotted part is $X^{s-1}$.\\
\begin{tikzpicture}[scale=0.5]


\draw[very thick](0,5)--(2,3); 
\draw[very thick](1,5)--(2,5);
\draw[very thick](1,5)--(0,4);
\draw[very thick](0,3)--(2,1);
\draw[very thick](0,2)--(2,4);
\draw[very thick](0,0)--(2,2);
\draw[very thick](0,1)--(1,0);
\draw[very thick](2,0)--(1,0);

\draw[very thick](2,6)--(2,-0.5);


\draw[very thick](2,5)--(4,3); 
\draw[very thick](3,5)--(4,5);
\draw[very thick](3,5)--(2,4);
\draw[very thick](2,3)--(4,1);
\draw[very thick](4,4)--(2,2);
\draw[very thick](2,0)--(4,2);
\draw[very thick](2,1)--(3,0);
\draw[very thick](4,0)--(3,0);

\draw[very thick](4,6)--(4,-0.5);

\draw[very thick](4,5)--(6,3); 
\draw[very thick](5,5)--(6,5);
\draw[very thick](5,5)--(4,4);
\draw[very thick](4,3)--(6,1);
\draw[very thick](4,2)--(6,4);
\draw[very thick](4,0)--(6,2);
\draw[very thick](4,1)--(5,0);
\draw[very thick](6,0)--(5,0);	
\draw[very thick](4,6)--(4,-0.5);
\draw[very thick](6,6)--(6,-0.5);

\draw[very thick](6,5)--(8,3); 
\draw[very thick](7,5)--(8,5);
\draw[very thick](7,5)--(6,4);
\draw[very thick](6,3)--(8,1);
\draw[very thick](6,2)--(8,4);
\draw[very thick](6,0)--(8,2);
\draw[very thick](6,1)--(7,0);
\draw[very thick](8,0)--(7,0);

\draw[very thick](8,6)--(8,-0.5);

\draw[very thick, dashed](10,5)--(12,3); 
\draw[very thick, dashed](11,5)--(12,5);
\draw[very thick, dashed](11,5)--(10,4);
\draw[very thick, dashed](10,3)--(12,1);
\draw[very thick, dashed](10,2)--(12,4);
\draw[very thick, dashed](10,0)--(12,2);
\draw[very thick, dashed](10,1)--(11,0);
\draw[very thick, dashed](12,0)--(11,0);

\draw[very thick](10,6)--(10,-0.5);
\draw[very thick](0,-1)--(10,-1);
\node at (5,-1.5) {$Z=X^{s-1}$};

\draw[very thick](8,5)--(10,3); 
\draw[very thick](9,5)--(10,5);
\draw[very thick](9,5)--(8,4);
\draw[very thick](8,3)--(10,1);
\draw[very thick](8,2)--(10,4);
\draw[very thick](8,0)--(10,2);
\draw[very thick](8,1)--(9,0);
\draw[very thick](10,0)--(9,0);
\node at (-1.4,5) {$\phi_{-1}(1)$};
\node at (-1.4,4) {$\phi_{-1}(2)$};
\node at (-1.4,3) {$\phi_{-1}(6)$};
\node at (-1.4,2) {$\phi_{-1}(3)$};
\node at (-1.4,1) {$\phi_{-1}(5)$};
\node at (-1.4,0) {$\phi_{-1}(4)$};

\node at (12.5,5) {$1$};
\node at (12.5,4) {$2$};
\node at (12.5,3) {$3$};
\node at (12.5,2) {$4$};
\node at (12.5,1) {$5$};
\node at (12.5,0) {$6$};

\node at (17,5) {$\pi_{-1}(1)=2$};
\node at (17,4) {$\pi_{-1}(2)=4$};
\node at (17,3) {$\pi_{-1}(3)=1$};
\node at (17,2) {$\pi_{-1}(4)=6$};
\node at (17,1) {$\pi_{-1}(5)=3$};
\node at (17,0) {$\pi_{-1}(6)=5$};
\node at (11,-1) {$X$};
\end{tikzpicture}
\subsubsection{The quasitoric case $\mathcal Q(s,ns+1)$}\label{quasitoric plus un}
\begin{minipage}{0.7\textwidth}
	We also denote $\pi_1$, $\phi_1$ the bijections associated to the building block $Y$ (see \S \ref{notations usual closure}). We see that $\pi_1(1)=s$ and for $i>1$, $\pi_1(i)=i-1$; thus for every $i$, $\pi_1(i)\equiv i-1\ \ (s)$. So for $1\leq k\leq s$, $\phi_1(k)=\pi_1^{k-1}(1)\equiv 2-k\ (s)$. \\
	Hence $\phi_1(1)=1$\ and if \ $k\geq 2$,\ $\phi_1(k)=s+2-k.$
\end{minipage}
\begin{minipage}{0.3\textwidth}\ \ 
	\begin{tikzpicture}
	\draw [very thick] (0.5,2) -- (1.5,0);
	\draw [very thick] (0.5,1.5) -- (1.5,2);
	\draw [very thick](0.5,1) -- (1.5, 1.5);
	\draw [very thick](0.5,0.5) -- (1.5, 1);
	\draw [very thick](0.5,0) -- (1.5, 0.5);
	\node at (0,2){$\phi_1(1)$};
	\node at (0,1.5){$\phi_1(5)$};
	\node at (0,1){$\phi_1(4)$};
	\node at (0,0.5){$\phi_1(3)$};
	\node at (0,0){$\phi_1(2)$};
	\end{tikzpicture}
\end{minipage}
\subsubsection{The quasitoric case $\mathcal Q(s,ns-1)$}\label{quasitoric moins un}
Since $Y^s$ is a pure braid, $\pi_{s-1}=\pi_1^{-1}$, hence $\pi_{s-1}(k)\equiv k+1$ for $k<s$ and 
$\pi_{s-1}(s)=1$. It follows that $\phi_{s-1}(k)=k$.
\subsection{Valid pairs of crossing points of rosette braids}
\subsubsection{Notations}
$\bullet$ In \S \ref{ou est le starting point} we picked the top left point as the starting point of the closed braid; we denote it $x$.\\ 
$\bullet$ Let  $C(i,j)$ be the set of all crossing points of the intervals $I_i$ and  $I_j$; equivalently of the strands $\phi(i)$ and $\phi(j)$.\\
$\bullet$ For $i,j,k,l$, we let $$S\big((i,j),(k,l)\big)=-\sum \epsilon_p(x)\epsilon_q(x)$$ where the sum goes over all interleaved crossing points $p\in C(i,j)$ and $q\in C(k,l)$. \\
$\bullet$
Let $u$ be the horizontal coordinate on the plane of the   braid shadow.

\subsubsection{Between $2$ strands}

\begin{lem}\label{pentes successives}
	Let $(p_k)$  be the crossing points of $\phi(i)$ and $\phi(j)$, ordered so that $u(p_k)<u(p_{k+1})$. For each $k$, $ \epsilon_{p_k}(x)\epsilon_{p_{k+1}}(x)=-1$. 
\end{lem}\begin{proof}
	The slopes of $\phi(i)$ have opposite signs at $p_k$ and $p_{k+1}$ (same for the $\phi(j)$). For example we can have\\ \begin{tikzpicture}
	\draw[very thick, -stealth](0,0)--(1,1); 
	\draw[very thick, -stealth](0,1)--(1,0);
	\node at (-0.4,1){$\phi(i)$};
	\node at (-0.4,0){$\phi(j)$};
	\node at (0.5,0){$p_k$};
	
	\draw[very thick, -stealth](3,0)--(4,1); 
	\draw[very thick, -stealth](3,1)--(4,0);
	\node at (2.6,1){$\phi(j)$};
	\node at (2.6,0){$\phi(i)$};
	\node at (3.5,-0.1){$p_{k+1}$};
	
	\end{tikzpicture}
\end{proof}
\begin{lem}\label{deux brins}
	Let $i<j$. 
	Then $S\big((i,j),(i,j)\big)=n$ thus\\ $\sum_{i\neq j}S\big((i,j),(i,j)\big)= \dfrac{s(s-1)}{2}n$.
	
\end{lem}

\begin{proof}
	Let $(p,\tilde{p})$, $(q,\tilde{q})$ two crossing points of the $i$-th and $j$-th intervals. If $u(p)<u(q)$, the crossing points $(p,\tilde{p})$, $(q,\tilde{q})$ are clearly interleaved.\\
	The cardinal number of $C(i,j)$ is equal to $2n$ or $2n+1$. Following Lemma \ref{pentes successives}, we can assume, without loss of generality, that the sequence of signs of the crossing points is 
	$$\underbrace{(+-)(+-)(+-)...(+-)(+-)}_\text{n times}\ \ \ \ \ \mbox{or}\ \ \ \ \ \underbrace{(+-)(+-)(+-)...(+-)(+-)}_\text{n times}+$$
	If $p$ is in a given block $(+-)$, the sum of the $\epsilon_p(x)\epsilon_q(x)$'s for the $q$'s in the following $(+-)$ blocks cancel out so we are left with  $\epsilon_p(x)\epsilon_q(x)$ for $p$ and the $q$ inside the same block and $\epsilon_p(x)\epsilon_q(x)=-1$ (Lemma \ref{pentes successives}).  Thus the total sum $S\big((i,j),(i,j)\big)$ is the number of blocks, i.e. $n$.
	
\end{proof}
\subsubsection{Between $3$ strands}
First we investigate when two crossing points involving exactly $3$ strands are interleaved.
\begin{lem}\label{le compte pour trois}
	Let $1\leq i<j<k\leq s$  and consider two crossing points $(p,\tilde{p})$, $(q,\tilde{q})$ which belong to two different sets among $C(i,j), C(i,k), C(j,k)$. Assume that $u(p)<u(q)$. The $(p,\tilde{p})$, $(q,\tilde{q})$ are interleaved  if and only if one of the following is true
	\begin{enumerate}
		\item $(p,\tilde{p})\in C(i,j)$,\ $(q,\tilde{q})\in C(i,k)$ 
		\item $(p,\tilde{p})\in C(j,k)$,\ $(q,\tilde{q})\in C(i,j)$ 
		\item $(p,\tilde{p})\in C(i,k)$,\ $(q,\tilde{q})\in C(j,k)$ 
	\end{enumerate}
	
\end{lem}
\begin{proof}
Two of the points in $\{p,\tilde{p},q,\tilde{q}\}$ belong to the same interval among $i,j,k$; the remaining two points belong to the remaining two intervals in $\{i,j,k\}$. If moreover $\big(p,\tilde{p}\big),\big(q,\tilde{q}\big)$ are interleaved, we have the following possibilities \\
\begin{minipage}{0.8\textwidth}
	\begin{tikzpicture}
	\node at (-1,0){$1.$};
	\draw [very thick](0,0) -- (3,0);
	\node at (0,0){$\bullet$};
	\node at (1,0){$\times$};
	\node at (4.5,0){$\times$};
	\node at (7.5,0){$\times$};
	\node at (2,0){$\times$};
	\node at (1,-0.4){$p$};
	\node at (2,-0.4){$q$};
	\node at (4.5,-0.4){$\tilde{p}$};
	\node at (7.5,-0.4){$\tilde{q}$};
	\node at (3,0){$\bullet$};
	\node at (1.5,0.5){i};
	\node at (4.5,0.5){j};
	\node at (7.5,0.5){k};
	\draw [very thick](6,0) -- (3,0);
	\node at (0,0){$\bullet$};
	\node at (6,0){$\bullet$};
	\draw [very thick](6,0) -- (9,0);
	\node at (0,0){$\bullet$};
	\node at (9,0){$\bullet$};
	
	\end{tikzpicture}
\end{minipage}
\begin{minipage}{0.2\textwidth}
	
\end{minipage}\\
\begin{minipage}{0.8\textwidth}
	\begin{tikzpicture}
	\node at (-1,0){$2.$};
	\draw [very thick](0,0) -- (3,0);
	\node at (0,0){$\bullet$};
	\node at (5,0){$\times$};
	\node at (1.5,0){$\times$};
	\node at (7.5,0){$\times$};
	\node at (7.5,-0.4){$\tilde{q}$};
	\node at (4,0){$\times$};
	\node at (1.5,-0.4){$p$};
	
	\node at (5,-0.4){$\tilde{p}$};
	\node at (4,-0.4){$q$};
	\node at (3,0){$\bullet$};
	\node at (1.5,0.5){i};
	\node at (4.5,0.5){j};
	\node at (7.5,0.5){k};
	\draw [very thick](6,0) -- (3,0);
	\node at (0,0){$\bullet$};
	\node at (6,0){$\bullet$};
	\draw [very thick](6,0) -- (9,0);
	\node at (0,0){$\bullet$};
	\node at (9,0){$\bullet$};
	
	\end{tikzpicture}
\end{minipage}
\begin{minipage}{0.2\textwidth}
	
\end{minipage}\\
\begin{minipage}{0.8\textwidth}
	\begin{tikzpicture}	
	\node at (-1,0){$3.$};
	\draw [very thick](0,0) -- (3,0);
	\node at (4.5,-0.4){$q$};
	\node at (4.5,0){$\times$};
	\node at (0,0){$\bullet$};
	\node at (1.5,0){$\times$};
	\node at (8,0){$\times$};
	\node at (8,-0.4){$\tilde{q}$};
	\node at (7,0){$\times$};
	\node at (7,-0.4){$\tilde{p}$};
	\node at (1.5,-0.4){$p$};
	\node at (3,0){$\bullet$};
	\node at (1.5,0.5){i};
	\node at (4.5,0.5){j};
	\node at (7.5,0.5){k};
	\draw [very thick](6,0) -- (3,0);
	\node at (0,0){$\bullet$};
	\node at (6,0){$\bullet$};
	\draw [very thick](6,0) -- (9,0);
	\node at (0,0){$\bullet$};
	\node at (9,0){$\bullet$};
	\end{tikzpicture}
\end{minipage}
\begin{minipage}{0.2\textwidth}
	
\end{minipage}
\end{proof}
\begin{notat}
		If $a\neq b$, $(ab)$ denotes a crossing point between the interval $I_a$ and $I_b$.
		\end{notat}
\begin{lem}\label{les deux types}
	Let $(i,j,k)$, $i<j<k$. For both the rosette and quasitoric braids, for $(s,ns+1)$ (resp. $s,ns-1)$ the crossing points are in the order
	\begin{equation}\label{ordre des points}
	...\underset{+}{(jk)}\underset{+}{(ik)}\underset{+}{(ij)}\underset{-}{(jk)}\underset{-}{(ik)}\underset{-}{(ij)}...(\mbox{resp.}...\underset{+}{(ij)}\underset{+}{(ik)}\underset{+}{(jk)}\underset{-}{(ij)}\underset{-}{(ik)}\underset{-}{(jk)}...)
	\end{equation}
	Note: the first crossing point of the braid can be any one among $\{(ij),(ij),(jk)\}$, this is why we wrote $...$ in (\ref{ordre des points}).\\ Under a crossing point $(ab)$ we wrote its sign $\epsilon$.  
\end{lem}
\begin{proof} I) {\bf The rosette case}\\
We use Lemma \ref{forme de u}.	We recall  that the $a$-th strand is going upwards (resp. downwards)  if $a$ is even (resp. odd).\\
	a)  $\mathbf{\mathcal R(s,ns+1)}$\\
	We go through all the cases.\\
	\\
\begin{minipage}{0.7\textwidth} 1) $i<j<k\leq A$, then $\phi(i), \phi(j)$ and $\phi(k)$ are all odd and $\phi(i)< \phi(j)<\phi(k)$\\
	$\underset{-}{(jk)}\underset{-}{(ik)}\underset{-}{(ij)}$
	\end{minipage}
\begin{minipage}{0.3\textwidth}
\ \ 	\begin{tikzpicture}
\draw [very thick] (1.5,0.5) -- (0,2);
\draw [very thick](0,1.5) -- (1, 0.5);
\draw [very thick](1.5,1) -- (1, 0.5);

\draw [very thick](0,1) -- (0.5,0.5);
\draw [very thick](1.5,1.5) -- (0.5,0.5);

\node at (-0.4,2){$\phi(i)$};
\node at (-0.4,1.5){$\phi(j)$};
\node at (-0.4,1){$\phi(k)$};
\end{tikzpicture}
\end{minipage}\\
\\
2) $i<j\leq A<k$. Then $\phi(k)$ is even, while $\phi(i)$ and $\phi(j)$ are odd with $\phi(i)<\phi(j)$. There are $3$ possibilities.\\
\begin{minipage}{0.7\textwidth} a) $\phi(k)<\phi(i)<\phi(j)$\\
	$\underset{-}{(ij)}\underset{+}{(jk)}\underset{+}{(ik)}$
\end{minipage}
\begin{minipage}{0.3\textwidth}
	\ \ 	\begin{tikzpicture}
	\draw [very thick] (0.5,2.5) -- (0,2);
		\draw [very thick] (0.5,2.5) -- (2.5,0.5);
	\draw [very thick](0,1.5) -- (1, 0.5);
	\draw [very thick](2.5,2) -- (1, 0.5);
	
	\draw [very thick](0,1) -- (0.5,0.5);
	\draw [very thick](2.5,2.5) -- (0.5,0.5);

	\node at (-0.4,2){$\phi(k)$};
	\node at (-0.4,1.5){$\phi(i)$};
	\node at (-0.4,1){$\phi(j)$};
	\end{tikzpicture}
\end{minipage}\\
\\
\\
\begin{minipage}{0.7\textwidth} b) $\phi(i)<\phi(k)<\phi(j)$\\
	$\underset{-}{(ik)}\underset{-}{(ij)}\underset{+}{(jk)}$
\end{minipage}
\begin{minipage}{0.3\textwidth}
	\ \ 	\begin{tikzpicture}
	\draw [very thick] (1.5,0.5) -- (0,2);
	
	\draw [very thick](0,1.5) -- (0.5, 2);
	\draw [very thick](2,0.5) -- (0.5, 2);

	\draw [very thick](0,1) -- (0.5,0.5);
	\draw [very thick](2,2) -- (0.5,0.5);

	\node at (-0.4,2){$\phi(i)$};
	\node at (-0.4,1.5){$\phi(k)$};
	\node at (-0.4,1){$\phi(j)$};
	\end{tikzpicture}
	\end{minipage}
	\\
	\\
	\begin{minipage}{0.7\textwidth} c) $\phi(i)<\phi(j)<\phi(k)$\\
		$\underset{-}{(jk)}\underset{-}{(ik)}\underset{-}{(ij)}$
	\end{minipage}
	\begin{minipage}{0.3\textwidth}
		\ \ 	\begin{tikzpicture}
		\draw [very thick] (1.5,0.5) -- (0,2);
		
	\draw [very thick](0,1.5) -- (1, 0.5);
		\draw [very thick](2,1.5) -- (1, 0.5);

		\draw [very thick](0,1) -- (1.5,2.5);

		\node at (-0.4,2){$\phi(i)$};
		\node at (-0.4,1.5){$\phi(j)$};
		\node at (-0.4,1){$\phi(k)$};
		\end{tikzpicture}
\end{minipage}\\
\\
3) $i\leq A<j<k$. Then $\phi(i)$ is odd, while $\phi(j)$ and $\phi(k)$ are even with $\phi(j)>\phi(k)$. There are $3$ possibilities.\\
\\
\begin{minipage}{0.7\textwidth} a) $\phi(i)<\phi(k)<\phi(j)$\\
	$\underset{-}{(ik)}\underset{-}{(ij)}\underset{+}{(jk)}$
\end{minipage}
\begin{minipage}{0.3\textwidth}
	\ \ 	\begin{tikzpicture}
	\draw [very thick] (1.5,0.5) -- (0,2);
	
	\draw [very thick](0,1.5) -- (1, 2.5);

	\draw [very thick](0,1) -- (1.5,2.5);
	
	\draw [very thick](2,1.5) -- (1,2.5);
	
	\node at (-0.4,2){$\phi(i)$};
	\node at (-0.4,1.5){$\phi(k)$};
	\node at (-0.4,1){$\phi(j)$};
	\end{tikzpicture}
\end{minipage}\\
\\
\\
\begin{minipage}{0.7\textwidth} b) $\phi(k)<\phi(i)<\phi(j)$\\
	$\underset{-}{(ij)}\underset{+}{(jk)}\underset{+}{(ik)}$
\end{minipage}
\begin{minipage}{0.3\textwidth}
	\ \ 	\begin{tikzpicture}
	\draw [very thick] (0.5,2.5) -- (0,2);
		\draw [very thick] (0.5,2.5) -- (2.5,0.5);

	\draw [very thick](2,1.5) -- (1,0.5);
	\draw [very thick](0,1.5) -- (1,0.5);
	
	\draw [very thick](0,1) -- (1.5,2.5);

	\node at (-0.4,2){$\phi(k)$};
	\node at (-0.4,1.5){$\phi(i)$};
	\node at (-0.4,1){$\phi(j)$};
	\end{tikzpicture}
\end{minipage}
\\
\\
\begin{minipage}{0.7\textwidth} c) $\phi(k)<\phi(j)<\phi(i)$\\
	$\underset{+}{(jk)}\underset{+}{(ik)}\underset{+}{(ij)}$
\end{minipage}
\begin{minipage}{0.3\textwidth}
	\ \ 	\begin{tikzpicture}

	\draw [very thick](0,1.5) -- (1, 2.5);
\draw [very thick](0,2) -- (0.5, 2.5);
\draw [very thick](2.5,1) -- (1, 2.5);

\draw [very thick](2.5,0.5) -- (0.5, 2.5);
	
	\draw [very thick](0,1) -- (0.5,0.5);

	\draw [very thick](2.5,2.5) -- (0.5,0.5);
	
	\node at (-0.4,2){$\phi(k)$};
	\node at (-0.4,1.5){$\phi(j)$};
	\node at (-0.4,1){$\phi(i)$};
	\end{tikzpicture}
\end{minipage}\\
\\
\begin{minipage}{0.7\textwidth}
4) $A<i<j<k$. Then  $\phi(i), \phi(j)$ and $\phi(k)$ are all even and $\phi(i)> \phi(j)>\phi(k)$\\
$\underset{+}{(jk)}\underset{+}{(ik)}\underset{+}{(ij)}$
	\end{minipage}
\begin{minipage}{0.3\textwidth}
		\begin{tikzpicture}
	\draw [very thick] (0.5,2.5) -- (0,2);
		\draw [very thick] (0.5,2.5) -- (2,1);
	\draw [very thick](0,1.5) -- (1, 2.5);
	
	\draw [very thick](2,1.5) -- (1, 2.5);

	\draw [very thick](0,1) -- (1.5,2.5);

	\node at (-0.4,2){$\phi(k)$};
	\node at (-0.4,1.5){$\phi(j)$};
	\node at (-0.4,1){$\phi(i)$};
	\end{tikzpicture}
	\end{minipage}
\\
\\
b)  $\mathbf{\mathcal R(s,ns-1)}$\\
Let $2\leq i<j<k$; we set $\tilde{k}=s+2-i$,  $\tilde{j}=s+2-j$ and $\tilde{i}=s+2-k$. Then $\tilde{i}<\tilde{j}<\tilde{k}$ and $\phi_{-1}(i)=\phi_1(\tilde{k})$, $\phi_{-1}(j)=\phi_1(\tilde{j})$ and $\phi_{-1}(k)=\phi_1(\tilde{i})$, cf. \S \ref{s,s moins un}.\\ We apply the  result above for $\mathcal R(s,sn+1)$: the $3$ intervals $\tilde{i},\tilde{j},\tilde{k}$ of $\mathcal R(s,ns+1)$ have crossing points in the following order \\ $$...\underset{+}{(\tilde{j}\tilde{k})}\underset{+}{(\tilde{i}\tilde{k})}\underset{+}{(\tilde{i}\tilde{j})}\underset{-}{(\tilde{j}\tilde{k})}\underset{-}{(\tilde{i}\tilde{k})}\underset{-}{(\tilde{i}\tilde{j})}\ \ \ \ \ (\star)$$. Thus the intervals $i,j,k$ of $\mathcal R(s,ns-1)$ have crossing points in the following order: 
	$$...\underset{-}{(ij)}\underset{-}{(ik)}\underset{-}{(jk)}\underset{+}{(ij)}\underset{+}{(ik)}\underset{+}{(jk)}\ \ \ \ \ (\star\star).$$ Indeed, since $a\mapsto s+2-a$ reverses the order of two integers, the signs $\epsilon(x)$ under the crossing points of $(\star)$ and $(\star\star)$ are opposite: for example if we have 
	\begin{minipage}{0.2\textwidth}
	for $\mathcal R(s,ns+1)$
	\end{minipage}
		\begin{minipage}{0.2\textwidth}
	\begin{tikzpicture}
	\draw[very thick, -stealth](0,0)--(1,1); 
	\draw[very thick, -stealth](0,1)--(1,0);
	\node at (-0.3,1){$\tilde{j}$};
	\node at (-0.3,0){$\tilde{k}$};
		\node at (0.5,-0.3){$\epsilon=-1$};
	\end{tikzpicture}
\end{minipage}
\begin{minipage}{0.4\textwidth}
	we have, for $\mathcal R(s,ns-1)$,
		\end{minipage}
\begin{minipage}{0.2\textwidth}
	\begin{tikzpicture}
	\draw[very thick, -stealth](3,0)--(4,1); 
	\draw[very thick, -stealth](3,1)--(4,0);
	\node at (2.7,1){$j$};
	\node at (2.7,0){$i$};
		\node at (3.5,-0.3){$\epsilon=+1$};
	
	\end{tikzpicture}
	\end{minipage}

	If $i=1$,   we proceed as above. We set $\tilde{i}=i=1$, $\tilde{k}=j$ and $\tilde{j}=k$
	and derive from $(\star)$
	 $$...\underset{-}{(jk)}\underset{+}{(ij)}\underset{+}{(ik)}\underset{+}{(jk)}\underset{-}{(ij)}\underset{-}{(ik)}...\ \ \ \ \ (\star\star\star).$$
	 The signs of the crossing points involving $i$ and $\tilde{i}$ remain the same and the sign of $(jk)$ is changed into its opposite.\\
	 \\
	 {\bf II) The quasitoric case}\\
	 \begin{minipage}{0.6\textwidth}
	 We illustrate the situation with $s=5$ and draw $X^5$. 
	 \end{minipage}
	 \begin{minipage}{0.4\textwidth}
	\begin{tikzpicture}
\draw [very thick] (0.5,2) -- (1.5,0);
\draw [very thick] (0.5,1.5) -- (1.5,2);
\draw [very thick](0.5,1) -- (1.5, 1.5);
\draw [very thick](0.5,0.5) -- (1.5, 1);
\draw [very thick](0.5,0) -- (1.5, 0.5);

\draw [very thick] (1.5,2) -- (2.5,0);
\draw [very thick] (1.5,1.5) -- (2.5,2);
\draw [very thick](1.5,1) -- (2.5, 1.5);
\draw [very thick](1.5,0.5) -- (2.5, 1);
\draw [very thick](1.5,0) -- (2.5, 0.5);

\draw [very thick] (3.5,2) -- (4.5,0);
\draw [very thick] (3.5,1.5) -- (4.5,2);
\draw [very thick](3.5,1) -- (4.5, 1.5);
\draw [very thick](3.5,0.5) -- (4.5, 1);
\draw [very thick](3.5,0) -- (4.5, 0.5);

\draw [very thick] (2.5,2) -- (3.5,0);
\draw [very thick] (2.5,1.5) -- (3.5,2);
\draw [very thick](2.5,1) -- (3.5, 1.5);
\draw [very thick](2.5,0.5) -- (3.5, 1);
\draw [very thick](2.5,0) -- (3.5, 0.5);

\draw [very thick] (4.5,2) -- (5.5,0);
\draw [very thick] (4.5,1.5) -- (5.5,2);
\draw [very thick](4.5,1) -- (5.5, 1.5);
\draw [very thick](4.5,0.5) -- (5.5, 1);
\draw [very thick](4.5,0) -- (5.5, 0.5);

\node at (0.2,0.5){$b$};
\node at (0.2,0){$c$};
\node at (0.2,1.5){$a$};

\end{tikzpicture}
 \end{minipage}
We consider $3$ strands (not intervals) $a<b<c$. The strand $a$ is the first one of the $3$ to reach the top point thus the first crossing point of strands in $\{a,b,c\}$ involves $a$ and $a$ meets first $b$ then $c$. Thus the order is
\begin{equation}\label{strands}
(ab)(ac)(bc)
\end{equation} 
\begin{minipage}{0.3\textwidth}
where 
\end{minipage}
\begin{minipage}{0.7\textwidth}

\begin{tikzpicture}
\draw[very thick, -stealth](0,0)--(1,1); 
\draw[very thick, -stealth](0,1)--(1,0);
\node at (-0.3,1){$a$};
\node at (-0.3,0){$b$};

\draw[very thick, -stealth](3,0)--(4,1); 
\draw[very thick, -stealth](3,1)--(4,0);
\node at (2.7,1){$a$};
\node at (2.7,0){$c$};

\draw[very thick, -stealth](6,0)--(7,1); 
\draw[very thick, -stealth](6,1)--(7,0);
\node at (5.7,1){$b$};
\node at (5.7,0){$c$};

\end{tikzpicture}

\end{minipage}
a) The case $(s,ns+1)$. \\

Consider $3$ intervals $2\leq i<j<k$: they correspond to strands $a<b<c$ with $i\mapsto c$, $j\mapsto b$ and $k\mapsto a$ (see \S \ref{quasitoric plus un}). Following (\ref{strands}), we have 
$\underset{+}{(jk)}\underset{+}{(ik)}\underset{+}{(ij)}$.\\
Suppose now that $i=1$, then $i\mapsto a$, $j\mapsto c$, $k\mapsto b$ and (\ref{strands}) yields $\underset{-}{(ik)}\underset{-}{(ij)}\underset{+}{(jk)}\underset{+}{(ik)}\underset{+}{(ij)}\underset{-}{(jk)}$, we recover $...\underset{+}{(jk)}\underset{+}{(ik)}\underset{+}{(ij)}...$.\\
b) The case $(s,ns-1)$. \\
If $i<j<k$, we have $i\mapsto a$, $j\mapsto b$, $k\mapsto c$ (see \S \ref{quasitoric moins un}) and (\ref{strands}) becomes $\underset{-}{(ij)}\underset{-}{(ik)}\underset{-}{(jk)}$. \end{proof}
\subsubsection{The count inside a block $X^s$ or $Y^s$}\label{count inside Xs}
\begin{itemize}
	\item {\bf For $(s,ns+1)$}. 
We do not know exactly which pair of crossing points is the first one so all we can derive from Lemma \ref{les deux types} is that the sequence of crossing points involving $i<j<k$ in a $X^s$ follows one of the $3$ patterns\\
$ \underset{a}{(jk)}\ \underset{a}{(ik)}\ \underset{a}{(ij)}\ \underset{-a}{(jk)}\ \underset{-a}{(ik)}\ \underset{-a}{(ij)},\ \ \ \ \ \  \  \ \underset{a}{(ik)}\ \underset{a}{(ij)}\ \underset{-a}{(jk)}\ \underset{-a}{(ik)}\ \underset{-a}{(ij)}\underset{a}{(jk)},\ \ \ \ \ \  \ \ \underset{a}{(ij)}\ \underset{-a}{(jk)}\ \underset{-a}{(ik)}\ \underset{-a}{(ij)}\underset{a}{(jk)}\underset{a}{(ik)}$
where $a\in\{-1,1\}$. Following Lemma \ref{le compte pour trois}, we go through the interleaved pairs of double points; for each of the $3$ patterns, we get\\. $-\dfrac{1}{4}\sum\epsilon_p(x)\epsilon_q(x)=+\dfrac{1}{4}$. 
\item {\bf For $(s,ns-1)$}.
We proceed similarly as above. The triples are \\ 
$ \underset{a}{(ij)}\ \underset{a}{(ik)}\ \underset{a}{(jk)}\ \underset{-a}{(ij)}\ \underset{-a}{(ik)}\ \underset{-a}{(jk)},\ \ \ \ \ \  \  \ \underset{a}{(ik)}\ \underset{a}{(jk)}\ \underset{-a}{(ij)}\ \underset{-a}{(ik)}\ \underset{-a}{(jk)}\underset{a}{(ij)},\ \ \ \ \ \  \ \ \underset{a}{(jk)}\ \underset{-a}{(ij)}\ \underset{-a}{(ik)}\ \underset{-a}{(jk)}\underset{a}{(ij)}\underset{a}{(ik)}$)
\\
We apply Lemma \ref{le compte pour trois} and get $-\dfrac{1}{4}\sum\epsilon_p(x)\epsilon_q(x)=-\dfrac{1}{4}$.	
	\end{itemize}
\subsubsection{The last piece}
We check that there is no contribution from the last piece of the braid, which is $X$ or $Y$ (resp. $X^{s-1}$ or $Y^{s-1}$) for $\mathcal R(s,sn+1)$ or $\mathcal Q(s,sn+1)$ (resp. $\mathcal R(s,sn-1)$ or $\mathcal Q(s,sn-1)$). 
\begin{itemize}
	\item {\bf For $\mathcal R(s,ns+1)$} (see Lemma \ref{forme de u}).
	Here are the two possibilities to have $2$ crossing points involving $3$ intervals $i,j,k$ in $X$. Note that in the figures below, the letters at the left of the strands indicate the indices of these strands.\\

\begin{minipage}{0.7\textwidth}
		{\bf 1) $\mathbf a$ is odd, $\mathbf{a+1}$, $\mathbf{a+3}$ are even}.\\ Then $a=\phi_1(i)$, $a+3=\phi_1(j)$ and $a+1=\phi_1(k)$. \\
		Thus the crossing points are in the order $(ik)(ij)$, which does not appear in the list of Lemma \ref{le compte pour trois}. 
	\end{minipage}
\begin{minipage}{0.3\textwidth}
	\ 	
		\begin{tikzpicture} 
		\draw [very thick] (1,1) -- (0,2);

		\draw [very thick](0.5,1) -- (1,1.5);
		
		\draw [very thick] (0,0.5) -- (1,1.5);
		\draw [very thick] (0,1.5) -- (0.5,2);
		
		\node at (-0.6,2){a};
		\node at (-0.6,1.5){a+1};
		\node at (-0.6,0.5){a+3};
		\end{tikzpicture}

	\end{minipage}\\

	\begin{minipage}{0.7\textwidth}
			{\bf 1) $\mathbf a$ is even, $\mathbf{a-1}$, $\mathbf{a-3}$ are odd}.\\
	Then $a=\phi_1(k), a-1=\phi_1(j), a-3=\phi_1(i)$. The crossing points are in the order $(jk)(ik)$, which does not appear in the list of Lemma \ref{le compte pour trois}. .
	\end{minipage}
\begin{minipage}{0.3\textwidth}
\ \ 		\begin{tikzpicture} 
	\draw [very thick] (1,1) -- (0,2);

	\draw [very thick](0.5,1) -- (1,1.5);
	
	\draw [very thick] (0,0.5) -- (1,1.5);
	\draw [very thick] (0,1) -- (0.5,0.5);
	
	\node at (-0.6,2){a-3};
	\node at (-0.6,1){a-1};
	\node at (-0.6,0.5){a};
	\end{tikzpicture}

\end{minipage}

			\item {\bf For $\mathcal Q(s,ns+1)$}.
			The only ordered pair of crossing points in $Y$ is of the form $(1k)(1j)$, i.e. $(ik)(ij)$, but $(ik)(ij)$ does not appear in the list of Lemma \ref{le compte pour trois}.
	
	\item {\bf For $\mathcal R(ns,s-1)$ or $\mathcal Q(ns,s-1)$}.\\
	Both $X^s$ and $Y^s$ contains $6$ crossing points involving $i,j,k$. 
	At most $2$ of these crossing points are contained in the last piece $X$ or $Y$. Thus  $X^{s-1}$ and $Y^{s-1}$ contains $4$,$5$ or $6$ crossing points involving $i,j,k$.\\
	If it contains  $6$ points, we have as above $-\dfrac{1}{4}\sum\epsilon_p(x)\epsilon_q(x)=-\dfrac{1}{4}$. \\
	We use Lemma \ref{les deux types} to list all the possible configurations if $X^{s-1}$ or $Y^{s-1}$ containing $4$ or $5$ crossing points for $i,j,k$. \\ 
	$\underset{+}{(ij)}\underset{+}{(ik)}\underset{+}{(jk)}\underset{-}{(ij)}\underset{-}{(ik)}$ \\
	$\underset{+}{(ik)}\underset{+}{(jk)}\underset{-}{(ij)}\underset{-}{(ik)}\underset{-}{(jk)}$ \\
	$\underset{-}{(jk)}\underset{+}{(ij)}\underset{+}{(ik)}\underset{+}{(jk)}\underset{-}{(ij)}$\\
	$\underset{+}{(ij)}\underset{+}{(ik)}\underset{+}{(jk)}\underset{-}{(ij)}$ \\
	$\underset{+}{(ik)}\underset{+}{(jk)}\underset{-}{(ij)}\underset{-}{(ik)}$ \\
	$\underset{-}{(jk)}\underset{+}{(ij)}\underset{+}{(ik)}\underset{+}{(jk)}$
	
	\end{itemize}
	For each of these cases, we sum over all the pairs given us by Lemma \ref{le compte pour trois} and we get every time $-\dfrac{1}{4}\sum\epsilon_p(x)\epsilon_q(x)=-\dfrac{1}{4}$.
	
	\subsection{Between $4$ strands}
	\begin{lem}
		Let $i<j$, $k<l$ and $\# \{i,j,k,l\}=4$.\\
1)  Two crossing points $(p,\tilde{p})\in\mathcal C(i,j)$ and $(q,\tilde{q})\in\mathcal C(k,l)$ are interleaved if and only if
			\begin{equation}\label{4 alternes}
			i<k<j<l.
			\end{equation}
2) In all $4$ cases, $\mathcal R(s,ns+1)$, $\mathcal Q(s,ns+1)$, $\mathcal R(s,ns-1)$, $\mathcal Q(s,ns-1)$	\begin{equation}\label{quatre}
		\ \sum_{p\in C(i,j),q\in C(k,l) }\epsilon_p(x)\epsilon_q(x)=0			\end{equation}

	\end{lem}
	\begin{proof} The inequalities (\ref{4 alternes}) are clear.
		In that case 
		\begin{equation}\label{equ.quatre}
		\sum_{p\in C(i,j),q\in C(k,l) }\epsilon_p(x)\epsilon_q(x)=\underbrace{\big(\sum_{p\in C(i,j)}\epsilon_p(x)\big)}_{(\star)}\underbrace{\big(\sum_{q\in C(k,l)}\epsilon_q(x)\big)}_{(\star\star)}.
		\end{equation}
		Thus (\ref{equ.quatre}) is non zero if both $(\star)$ and $(\star\star)$ are non zero. 
		Note that $X^s$, $Y^s$ contain $2$ points in $C(i,j)$ of opposite signs and the same is true for $C(k,l)$. Hence, for $\mathcal R(sm+1)$, $\mathcal Q(sm+1)$ (resp. $\mathcal R(sm-1)$, $\mathcal Q(sm+1)$) (\ref{equ.quatre}) is non zero if $X$, $Y$ (resp. $X^{s-1}$, $Y^{s-1}$) contains a single point in $C(i,j)$ and a single point in $C(k,l)$. But $X^{s-1}$, $Y^{s-1}$ contains a single point in a $C(i,j)$ if and only if the same is true for  $X,Y$. So we just need to ask: can $X$ or $Y$ contain a point in $C(i,j)$ and a point in $C(k,l)$?\\
		$\bullet$ {\bf The case of $X$}. Every crossing point in $X$ is between an odd numbered strand and an even numbered strand. Thus we derive from Lemma \ref{forme de u} that, for some $a,b$\\
		i) 
		$\phi_1(i)=2a-1$, $\phi_1(j)\in\{2a,2a+2\}$\\
		ii) $\phi_1(k)=2b-1$, $\phi_1(l)\in\{2b,2b+2\}$\\
		Since $i<k$, we have $a<b$. On the other hand, $j<l$ so $2b<2a+2$, i.e. $b\leq a$. We have a contradiction.
		\\
		$\bullet$ {\bf The case of $Y$}.
		It is clear that all the crossing points in $Y$ have a common strand.
			\end{proof}
	
\section{Another model for $3$-braids: proof of Proposition \ref{autre modele}}\label{3standard}
We require the braid shadow $\Gamma$ to close into a knot: for $3$-braids, this means that the 
permutation $\pi(\Gamma)$ induced by $\Gamma$ must have signature $1$ without being the identity. In particular $\Gamma$ is the product of an even number of braid shadow generators.
\begin{notat} Given a positive integer $M$, we let $B(3,2M)$ be the set of braid diagrams of length $2M$; we wiew such a diagram as a sequence of $M$ spots, each spot being filled by one of the following: $\hat{\sigma_1}^{2}$, $\hat{\sigma_2}^{2}$, $\hat{\sigma_2}\hat{\sigma_1}$, $\hat{\sigma_1}\hat{\sigma_2}$. For $\Gamma\in B(3,2M)$, we let $\alpha(\Gamma), \beta(\Gamma), \gamma(\Gamma), \delta(\Gamma)$ be the number of spots repectively filled  by $\hat{\sigma_1}^{2}$, $\hat{\sigma_2}^{2}$, $\hat{\sigma_2}\hat{\sigma_1}$, $\hat{\sigma_1}\hat{\sigma_2}$.\\
	Finally we let $KB(3,2M)$ the set of braid shadows in $B(3,2M)$ which close in a curve with a single component. 
	\end{notat}
Ma ([Ma]) proved that a random $n$-braid closes in a knot with probability $\dfrac{1}{n}$. In the case of $3$ braids with an even number of generators, this means a probability $\dfrac{2}{3}$. Thus $\#  B(3,2M)=4^M$ and $\#  KB(3,2M)=\dfrac{2}{3}4^M$.\\
\begin{lem}\label{debut du calcul}
For	$\Gamma\in KB(3,2M)$, $|defect(\Gamma)-2M|\leq 2|\gamma(\Gamma)-\delta(\Gamma)|$ 
\end{lem}
\begin{proof} 
	Chang and Erickson ([C-E1]) define three {\it homotopy moves} for the defect similarly to the Reidmeister moves. We recall the $2\longrightarrow 0$ move.\\
	\begin{tikzpicture}[scale=0.7]
	\draw 	(-0.3,-0.5) node {$p$};
	
	\draw [very thick](0,-1) arc (170:10:2);
	\draw [very thick](0,0) arc (10:170:-2);
	
	\draw [very thick](8,-2.5) arc (170:10:2);
	\draw [very thick](8,1.5) arc (10:170:-2);
	\draw [ultra thick] (5,-0.5)--(7,-0.5);
	\draw [ultra thick] (5,-0.5)--(7,-0.5);
	\draw [ultra thick] (6.5,-0.3)--(7,-0.5);
	\draw [ultra thick] (6.5,-0.7)--(7,-0.5);
	\draw 	(2,-2.2) node {$\textcircled{1}$};
	\draw 	(4.1,-0.5) node {$q$};
	\draw 	(10,-2) node {$\textcircled{2}$};
	
	\end{tikzpicture}\\
	If $p,q$ are interleaved, then $defect(\textcircled{1})=defect(\textcircled{2})+2$. \\
	For a $\hat{\sigma}_1^2$ or a $\hat{\sigma}_2^2$ in $\Gamma$, the corresponding two crossing points are interleaved hence we do a $2\longrightarrow 0$ move. We do these moves for all $\hat{\sigma}_1^2$'s, $\hat{\sigma}_2^2$'s and get a braid diagram $\Gamma_1$ with
	\begin{equation}\label{defect1}
	defect(\Gamma)=defect(\Gamma_1)+2\big(\alpha(\Gamma)+\beta(\Gamma)\big)
	\end{equation}
	The diagram $\Gamma_1$ is made of $\gamma(\Gamma_1)$ copies of $\hat{\sigma}_1\hat{\sigma}_2$ and $\delta(\Gamma_1)$ copies of $\hat{\sigma}_1\hat{\sigma}_2$.  
	If $\Gamma_1$ contains a $(\hat{\sigma}_2\hat{\sigma}_1)(\hat{\sigma}_1\hat{\sigma}_2)$, we do a $2\longrightarrow 0$ homotopy move which desingularizes  $\hat{\sigma}_1\hat{\sigma}_1$ and replaces $(\hat{\sigma}_2\hat{\sigma}_1)(\hat{\sigma}_1\hat{\sigma}_2)$ by $\hat{\sigma}_2\hat{\sigma}_2$ to which we apply another $2\longrightarrow 0$ move. We proceed similarly with the $(\hat{\sigma}_1\hat{\sigma}_2)(\hat{\sigma}_2\hat{\sigma}_1)$'s contained in $\Gamma_1$.\\
	After all these moves, we end up with a braid shadow $\Gamma_2$ made up entirely of 
	$\hat{\sigma}_2\hat{\sigma}_1$'s or entirely of $\hat{\sigma}_1\hat{\sigma}_2$'s and of length $2|\gamma(\Gamma)-\delta(\Gamma)|$. We have
	\begin{equation}\label{defect2}
	defect(\Gamma_1)=defect(\Gamma_2)+2(\gamma(\Gamma)+\delta(\Gamma)-|\gamma(\Gamma)-\delta(\Gamma)|)
	\end{equation}
	Note that the homotopy moves above did not change the permutation $\pi(\Gamma)$. Thus the curve $\Gamma_2$ is a $(3,2|\gamma(\Gamma)-\delta(\Gamma)|)$ quasitoric/rosette braid which closes in a knot; we derive from  Theorem \ref{thm fermeture usuelle} that
	\begin{equation}\label{defect3} |defect(\Gamma_2)|\leq 2|\gamma(\Gamma)-\delta(\Gamma)| 
	\end{equation}
	Putting together (\ref{defect1}), (\ref{defect2}) and (\ref{defect3}) yield the Lemma.
\end{proof} 
We now sum the $|\gamma(\Gamma)-\delta(\Gamma)| $'s over $\Gamma\in B(3,2M)$. We start by counting the  $\Gamma$'s corresponding to the same number $|\gamma(\Gamma)-\delta(\Gamma)|$.\\
1) we pick $k=\gamma(\Gamma)+\delta(\Gamma)\in\{0,2M\}$\\
2) given $k$, we pick the set $E_k$ of spots of $\hat{\sigma}_1\hat{\sigma}_2$, $\hat{\sigma}_2\hat{\sigma}_1$ in $\Gamma$; we have $\begin{pmatrix} M\\k\end{pmatrix}$ possibilities for $E_k$.\\
3) We can fill each spot in $\Gamma$ outside of $E_k$ by $\hat{\sigma}_1^2$ or $\hat{\sigma}_2^2$; this gives us a total of $2^{M-k}$ to fill the spots outside of $E_k$. \\
4) We count the different ways in which we can write $E_k$ as a product of  $\hat{\sigma_1}\hat{\sigma_2}$'s and $\hat{\sigma_2}\hat{\sigma_1}$'s. First we choose the number $\gamma=\gamma(\Gamma)\leq k$; for each such $\gamma$, we choose the set $A$ of spots for $\hat{\sigma_1}\hat{\sigma_2}$'s inside $E_k$ and this gives us $\begin{pmatrix}k\\ \gamma\end{pmatrix}$ possibilities.\\
\\
For each of the $\Gamma$'s that we just listed above, we add $|\gamma(\Gamma)-\delta(\Gamma)|=|k-2\gamma(\Gamma)|$ where $k$ was defined in 1). Going through 1)-4) above we end up with
\begin{equation}\label{grosse somme}
\sum_{\Gamma\in B(3,M)}|\gamma(\Gamma)-\delta(\Gamma)|=\sum_{k=1}^M\begin{pmatrix}M\\ k\end{pmatrix}2^{M-k}\underbrace{\sum_{\gamma=0}^k\begin{pmatrix}k\\  \gamma\end{pmatrix}|k-2\gamma|}_{(\star\star\star)}
\end{equation}
\begin{lem}\label{petits calculs}
$(\star\star\star)\leq 4\sqrt{k}2^k$
\end{lem}
\begin{proof} Since $\begin{pmatrix}k\\ \gamma\end{pmatrix}=\begin{pmatrix}k\\ k-\gamma\end{pmatrix}$ and $|k-2\gamma|=|k-2(k-\gamma)|$, we have
	 $$(\star\star\star)\leq 2\sum_{\gamma=0}^{[\frac{k}{2}]}\begin{pmatrix}k\\ \gamma\end{pmatrix}|k-2\gamma|.$$
	We recall identities from ([G-R]).\\
	a)	 $\sum_{a=0}^ma\begin{pmatrix}2m\\ a\end{pmatrix}=\dfrac{m}{2}2^{2m}$
	\\
	b) $\sum_{a=0}^{m}\begin{pmatrix}2m\\ a\end{pmatrix}=2^{2m-1}+\dfrac{1}{2}\begin{pmatrix}2m\\ m\end{pmatrix}$. \\
	If $k=2m$, a) and b)  yield $\sum_{a=0}^m(2m-2a)\begin{pmatrix}2m\\\ a\end{pmatrix}=m\begin{pmatrix}2m\\ m\end{pmatrix}$.\\
	If $k=2m+1$,\\ $\sum_{a=0}^m(2m+1-2a)\begin{pmatrix}2m+1\\\ a\end{pmatrix}\leq \sum_{a=0}^{m+1}(2m+2-2a)\begin{pmatrix}2m+2\\\ a\end{pmatrix}$ 
	$$\leq (m+1)\begin{pmatrix}2m+2\\ m+1\end{pmatrix}.$$
	
	We conclude using Sasv\'{a}ri's inequality 
	\begin{thm}([Sa])\label{sasvari}
		$	\begin{pmatrix}2s
		\\\ s\end{pmatrix}<\dfrac{4^s}{\sqrt{s\pi}}exp(-\dfrac{1}{8s}+\dfrac{1}{192s^3})<\dfrac{2^{2s}}{\sqrt{2s}}$
	\end{thm}
\end{proof}
We plug Lemma \ref{petits calculs} into (\ref{grosse somme}) and derive
\begin{equation}\label{suite du calcul}
\sum_{\Gamma\in B(3,2M)}|\gamma(\Gamma)-\delta(\Gamma)|\leq 4\sum_{k=1}^M\begin{pmatrix}M\\ k\end{pmatrix}2^{M}\sqrt{k}\leq 4^{M+1}\sqrt{M}
\end{equation}
We let $E(defect)$ be the expectation of the defect over $KB(3,2M)$; Lemma \ref{debut du calcul} and (\ref{suite du calcul}) yield
\begin{equation}\label{presque}
|E(defect)-2M|\leq \dfrac{4^{M+1}\sqrt{M}}{\#KB(3,2M)}=6\sqrt{M}
\end{equation}
We divide (\ref{presque}) by $16M=8c$ and compute the expectation of $v_2$ over the knots having diagrams in $KB(3,2M)$. We get $|\dfrac{E(v_2)}{c}-\dfrac{1}{8}|< \dfrac{1}{\sqrt{c}}$ which concludes the proof.
\section{Checkerboard diagram}
	\begin{tikzpicture}
	\draw (-1,4.5) node{$\bullet$};
		\draw (-1.3,4.5) node{$x$};
\draw 	(3.5,5.5) node {$0$};
\draw 	(3.5,4.5) node {$+1$};
\draw 	(3.5,2.5) node {$-1$};
\draw 	(3.5,0.5) node {$+1$};
\draw 	(-0.5,4.5) node {$+1$};
\draw 	(-0.5,2.5) node {$-1$};
\draw 	(-0.5,0.5) node {$+1$};

\draw 	(1.5,4.5) node {$-1$};
\draw 	(1.5,2.5) node {$+1$};
\draw 	(1.5,0.5) node {$-1$};
\draw 	(5.5,4.5) node {$-1$};
\draw 	(5.5,2.5) node {$+1$};
\draw 	(5.5,0.5) node {$-1$};
\draw 	(3.5,5.5) node {$0$};
\draw 	(7.5,4.5) node {$+1$};
\draw 	(7.5,2.5) node {$-1$};
\draw 	(7.5,0.5) node {$+1$};

\draw 	(0.5,3.5) node {$0$};
\draw 	(0.5,1.5) node {$0$};
\draw 	(2.5,3.5) node {$0$};
\draw 	(2.5,1.5) node {$0$};
\draw 	(4.5,3.5) node {$0$};
\draw 	(4.5,1.5) node {$0$};
\draw 	(6.5,3.5) node {$0$};
\draw 	(6.5,1.5) node {$0$};

\draw [very thick] (0,0) -- (5,5);
\draw [very thick][-stealth]  (0,0) -- (1,1);
\draw [very thick][-stealth]  (0,0) -- (3,3);
\draw [very thick] (-1,0) -- (0,0);
\draw [very thick] [-stealth] (-1,0) -- (-0.5,0);
\draw [very thick] (6,5) -- (5,5);
\draw [very thick] (5,5) -- (6,5);
\draw [very thick] (8,3) -- (6,5);
\draw [very thick][-stealth] (8,3) -- (8,2);
\draw [very thick][-stealth] -- (6,5)--(7,4);
\draw [very thick][-stealth] (5,5) -- (6,5);
\draw [very thick] (8,2) -- (6,0);
\draw [very thick][-stealth]  (8,2) -- (7,1);
\draw [very thick] [-stealth] (6,0) -- (5,0);
\draw [very thick]  (0,5) -- (5,0);
\draw [very thick] [-stealth] (5,0) -- (4,1);
\draw [very thick] [-stealth] (4,1)--(2.5,2.5);
\draw [very thick] [-stealth](2.5,2.5)--(1,4);
\draw [very thick] [-stealth] (0,5)--(-1,5);
\draw [very thick] [-stealth] (-1,5)--(-1,4);
\draw [very thick] [-stealth] (-1,4)--(0,3);
\draw [very thick] [-stealth] (0,3)--(2,1);
\draw [very thick] (3,0)--(2,1);
\draw [very thick] [-stealth] (3,0)--(4,0);
\draw [very thick] [-stealth] (4,0)--(8,4);
\draw [very thick] [-stealth] (4,0)--(6,2);
\draw [very thick] [-stealth] (4,0)--(7,3);
\draw [very thick] [-stealth] (8,4)--(8,5);
\draw [very thick] [-stealth] (8,5)--(7,5);
\draw [very thick] [-stealth] (7,5)--(2,0);
\draw [very thick] [-stealth] (7,5)--(6,4);
\draw [very thick] [-stealth] (6,4)--(5,3);
\draw [very thick] [-stealth] (6,4)--(3,1);
\draw [very thick] [-stealth] (2,0)--(1,0);
\draw [very thick] [-stealth] (1,0)--(-1,2);
\draw [very thick] [-stealth] (1,0)--(0,1);
\draw [very thick] [-stealth] (-1,2)--(-1,3);
\draw [very thick] [-stealth] (-1,3)--(1,5);
\draw [very thick] [-stealth] (-1,3)--(0,4);
\draw [very thick] [-stealth] (1,5)--(2,5);
\draw [very thick] [-stealth] (2,5)--(3,4);
\draw [very thick] [-stealth] (2,5)--(4,3);
\draw [very thick] [-stealth] (2,5)--(7,0);
\draw [very thick] [-stealth] (2,5)--(6,1);
\draw [very thick] [-stealth] (7,0)--(8,0);
\draw [very thick] [-stealth] (8,0)--(8,1);
\draw [very thick] [-stealth] (8,1)--(4,5);
\draw [very thick] [-stealth] (8,1)--(6,3);
\draw [very thick] [-stealth] (8,1)--(5,4);
\draw [very thick] [-stealth] (4,5)--(3,5);
\draw [very thick] [-stealth] (3,5)--(-1,1);
\draw [very thick] [-stealth] (3,5)--(2,4);
\draw [very thick] [-stealth] (3,5)--(1,3);
\draw [very thick] [-stealth] (-1,1)--(-1,0);
\draw [very thick](-0.4,0.6) circle (0.5);
\draw [very thick](-0.4,2.5) circle (0.5);
\draw [very thick](-0.4,0.6) circle (0.5);
\draw [very thick](-0.4,2.5) circle (0.5);
\draw [very thick](-0.4,4.4) circle (0.5);
\draw [very thick](-0.4,2.5) circle (0.5);
\draw [very thick](-0.4,2.5) circle (0.5);
\draw [very thick](-0.4,0.6) circle (0.5);
\draw [very thick](-0.4,2.5) circle (0.5);
\draw [very thick](1.5,4.4) circle (0.5);
\draw [very thick](1.5,2.5) circle (0.5);
\draw [very thick](1.5,0.6) circle (0.5);
\draw [very thick](3.5,4.4) circle (0.5);
\draw [very thick](3.5,2.5) circle (0.5);
\draw [very thick](3.5,0.6) circle (0.5);

\draw [very thick](5.5,4.4) circle (0.5);
\draw [very thick](5.5,2.5) circle (0.5);
\draw [very thick](5.5,0.6) circle (0.5);
\draw [very thick] [stealth-] (-0.5,1.1)--(-0.4,1.1);
\draw [very thick] [-stealth] (-0.5,3)--(-0.4,3);
\draw [very thick] [stealth-] (-0.5,4.9)--(-0.4,4.9);
\draw [very thick](7.4,4.4) circle (0.5);
\draw [very thick](7.4,2.5) circle (0.5);
\draw [very thick](7.4,0.6) circle (0.5);

\draw [very thick] [-stealth] (1.5,1.1)--(1.6,1.1);
\draw [very thick] [stealth-] (1.5,3)--(1.6,3);
\draw [very thick] [-stealth] (1.5,4.9)--(1.6,4.9);

\draw [very thick] [stealth-] (3.5,1.1)--(3.6,1.1);
\draw [very thick] [-stealth] (3.5,3)--(3.6,3);
\draw [very thick] [stealth-] (3.5,4.9)--(3.6,4.9);
\draw [very thick] [-stealth] (5.5,1.1)--(5.6,1.1);
\draw [very thick] [stealth-] (5.5,3)--(5.6,3);
\draw [very thick] [-stealth] (5.5,4.9)--(5.6,4.9);

\draw [very thick] [stealth-] (7.5,1.1)--(7.6,1.1);
\draw [very thick] [-stealth] (7.5,3)--(7.6,3);
\draw [very thick] [stealth-] (7.5,4.9)--(7.6,4.9);
\draw 	(4,-0.7) node {$p\ \mbox{Seifert circles}$};
\draw 	(-3,3) node {$q\ \mbox{Seifert circles}$};
\end{tikzpicture}\\
\\
We use the figure above to explain the proof which is based on Arnold's invariants (\S \ref{arnold invariants}). \\ The starting point is $x$, the oriented circles are the Seifert circles and the integers $\pm 1, 0$ denote the indices of the various regions.\\
We see that $ind_K(x)=\dfrac{1}{2}$; also, for every double point $d$, the average of the indices of the regions around $d$ is $0$. It follows that $St=0$.\\
We also see that for every Seifert circle $c$, $t(c)=0$ and that the total number of Seifert circles is $\boxed{s=pq}$.\\
We now count the number $\nu$ of crossing points. First consider a Seifert circle $c$ not on the last line and not on the last column and map it to  $2$ different crossing points, the one just below it and the one just on its right: this gives us $\boxed{(1)=2(pq-p-q+1)}$ crossing points. Next we let $c$ be one of the $p-1$ circles on the last line but not the last one on the right. We map it to the point immediately on the right of $c$ which gives us $\boxed{(2)=p-1}$ new crossing points. Finally we take for $c$ a Seifert circle which is on the right column but is not the bottom one and we map it to the crossing point just under it: this gives us $\boxed{(3)=q-1}$ crossing points. Adding these three quantities, we get 
$$\nu=(1)+(2)+(3)=2pq-p-q$$
Hence, $defect=J^+=1+\nu-s=(p-1)(q-1)$. \ \ \ \ \qed

Univ. Paris Est Creteil, CNRS, LAMA, F-94010 Creteil, France\\
villemarina@yahoo.fr
\end{document}